\documentclass{article}
\usepackage{amsmath,amssymb}
\newcommand{\dis}{\displaystyle}

\newcommand{\Div}{\mathrm{div}}
\newcommand{\dt}{{\partial_t}}

\newcommand{\R}{\mathbb{R}}

\def\qed{\hbox{${\vcenter{\vbox{
  \hrule height 0.4pt\hbox{\vrule width 0.4pt height 6pt
  \kern5pt\vrule width 0.4pt}\hrule height 0.4pt}}}$}}

\newtheorem{theo}{Theorem}

\newtheorem{lemm}[theo]{Lemma}

\begin{document}

\title{Regularity criterion for 3D Navier-Stokes equations in terms of the direction of the velocity\thanks{This work  was supported in part by NSF Grant
DMS-0607953.}}
%\titlerunning{Criterion for regularity to Navier-Stokes}

\author{Alexis  Vasseur}

%\institute{A. F. Vasseur \at
%               University of Texas at Austin, department of Mathematics,
%1 University Station C1200, Austin, TX78712-0257\\
%              Tel.: +512-471-2363\\
%              Fax: +512-471-9038\\
%              \email{vasseur@math.utexas.edu}   }



\maketitle

\bibliographystyle{plain}

\begin{abstract}
In this short note, we give  a link between the regularity of the
solution $u$ to the 3D Navier-Stokes equation, and  the behavior of
the direction of the velocity $u/|u|$. It is shown that the control
of $\Div (u/|u|)$ in a suitable $L_t^p(L_x^q)$ norm is enough to
ensure global regularity. The result is reminiscent of the criterion
in terms of the direction of the vorticity, introduced first by
Constantin and Fefferman. But in this case the condition is not on
the vorticity, but on the velocity itself. The proof, based on very
standard methods, relies on a straightforward relation between the
divergence of the direction of the velocity  and the growth of
energy along streamlines.
%\keywords{Navier-Stokes \and
%regularity criterion \and a priori estimates} \subclass{MSC 35B65
%\and 76D03 \and 76D05}
\end{abstract}

\section{Introduction}

This short paper deals with a new formulation  of well known
criterions for regularity of solutions to the incompressible
Navier-Stokes equation in dimension 3, namely:
\begin{equation}\label{eq_NS}
\begin{array}{l}
\dis{\dt u+\Div(u\otimes u)+\nabla P-\Delta u=0\qquad
t\in ]0,\infty[, \ x\in\R^3,}\\[3mm]
\dis{\mathrm{div} u=0.}
\end{array}
\end{equation}
The unknown is the velocity field $u(t,x)\in \R^3$. The pressure $P$
is a non local operator of $u$ which can be seen as a Lagrange
multiplier associated to the constraint of incompressibility $\Div
u=0$. The existence of weak solutions  was proved long ago by Leray
\cite {Leray} and Hopf \cite{Hopf}. They have  shown that, for any
initial value with finite energy $u^0\in L^2(\R^3)$,  there exists a
function $u\in L^\infty (0,\infty;L^2(\R^3))\times
L^2(0,\infty;\dot{H}^1(\R^3))$ verifying (\ref{eq_NS}) in the sense
of distribution, and verifying in addition the energy inequality:
\begin{equation}\label{eq_energy}
\|u(t,\cdot)\|^2_{L^2(\R^3)}+2\int_0^t\|\nabla
u(s,\cdot)\|^2_{L^2(\R^3)}\,ds\leq\|u_0\|^2_{L^2(\R^3)}, \qquad
t\geq0.
\end{equation}
Such a solution is now called Leray-Hopf weak solution to
(\ref{eq_NS}).
%From that time on, much effort has been made to
%establish results on the uniqueness and regularity of weak
%solutions. However those two questions remains yet mostly open.
%Especially it is not known until now if such a weak solution can
%develop singularities in finite time, even considering smooth
%initial data. The question of uniqueness is related to the one of
%regularity. Indeed it is well known that if the solution is smooth
%%enough, then  it is unique. Several steps has already been performed
%concerning the regularity of weak solutions.

In \cite{Serrin}, Serrin showed that a Leray-Hopf solution of
(\ref{eq_NS}) lying in $L^p(0,\infty;L^q(\R^3))$ with $p,q\geq1$
such that $2/p+3/q<1$ is smooth in the spatial directions. This
result was later extended in \cite{Struwe} and \cite{Fabes} to the
case of equality for $p<\infty$. Notice that the case of
$L^\infty(0,\infty;L^3(\R^3))$ was proven only very recently by
Escauriaza, Seregin and Sverak \cite{Seregin}.

An other class of regularity criterion was introduced by Beir\~ao da
Vaiga \cite{Vaiga} which involves the gradient of $u$. More
precisely, he showed that any Leray-Hopf solutions $u$ such that
$\nabla u$ lies in $L^p(L^q)$ with $2/p+3/q=2$, $3/2<q<\infty$, is
smooth. Beale-Kato-Majda \cite{Beale} dealt with the vorticity
$\omega=\mathrm{rot}\ u$ and proved regularity under the condition
$\omega\in L^1(L^\infty)$. This condition was later improved to
$L^1(BMO)$ by Kozono and Taniuchi \cite{Kozono}.

In \cite{Constantin}, Constantin and Fefferman introduced a
criterion involving the direction of the volicity $\omega/|\omega|$.
They have showed that under a Lipshitz-like regularity assumption on
$\omega/|\omega|$, the solution is smooth. (see \cite{Zhou} for
extension of this result).

Our result is of the same spirit but involves the direction of the
velocity itself instead of the vorticity.
\begin{theo}\label{th_rough_theorem}
Let $u$ be a Leray-Hopf solution to Navier-Stokes equations with
initial value $u_0\in L^2(\R^3)$. If $\Div(u/|u|)\in
L^p(0,\infty;L^q(\R^3))$ with:
$$
\frac{2}{p}+\frac{3}{q}\leq\frac{1}{2},\qquad q\geq 6,\qquad p\geq
4,
$$
Then $u$ is smooth on $(0,\infty)\times\R^3$.
\end{theo}
The result shows that it is enough to control the rate of change of
the direction of the velocity to get full regularity of the
solution. The main point of this modest paper is the following
straightforward equality coming from the incompressibility of the
flow:
\begin{equation}\label{eq_incompressibility}
|u|\Div(u/|u|)=-\frac{u}{|u|}\cdot\nabla|u|.
\end{equation}
This equality shows that, due to the incompressibility, the growth
of $|u|$ along the stream lines is linked to the divergence of the
direction of $u$. It means that to allow some increase of kinetic
energy $|u|^2$ along the streamlines, those streamlines need to be
bent, producing some divergence on the direction of the velocity.

This remark is the main point of this short note. The proof of the
theorem  then follows in a very standard way. It uses the fact that
the right-hand side term in (\ref{eq_incompressibility})
corresponds, up to the multiplication by a power of $|u|$, to the
flux of energy $u\cdot\nabla|u|^2$. Besides, It is also interesting
noticing that this term depends only on the symmetric part of the
gradient of $u$. Indeed it can be rewritten:
\begin{eqnarray*}
|u|\Div(u/|u|)&=&-\frac{u}{|u|}\cdot\nabla|u|=-\frac{u}{2|u|^2}\cdot\nabla|u|^2\\
&=&-\frac{u^T}{|u|^2}\cdot\nabla u\cdot u=-\frac{u^T}{|u|^2}\cdot
D(u)\cdot u.
\end{eqnarray*}
It was already known that if one component of the velocity is
bounded in a suitable space, then the solution is smooth (see Penel
and Pokorny \cite{Penel}, He \cite{He}, Zhou \cite{Zhou}, Chae and
Choe \cite{Chae}). Our result states that if  the direction of the
velocity does not change too drastically, the conclusion is still
true.

\section{Proof of Theorem \ref{th_rough_theorem}}

Let us first state a technical lemma:
%\begin{lemm}\label{lemm_alpha=1}
%If $||u|^\gamma\Div(u/|u|)|\in L^p(0,\infty;L^q(\R^3))$,
%$1\leq\overline{p}\leq\infty$, $1\leq\overline{q}\leq\infty$, with
%$0\leq\gamma\leq1$ and:
%$$
%\frac{2}{p}+\frac{3}{q}+\frac{3(1-\gamma)}{2}\leq 5,
%$$
%%then there exists $\overline{p}\geq 1, \overline{q}$ such that:
%$$
%||u|\Div(u/|u|)|\in
%L^{\overline{p}}(0,\infty;L^{\overline{q}}(\R^3)),
%$$
%with:
%$$
%\frac{2}{\overline{p}}+\frac{3}{\overline{q}}=5-\frac{2}{p}-\frac{3}{q}-\frac{3\gamma}{2}.
%$$
%If either $p\neq \infty$ or $\gamma\neq1$, then we can choose
%$\overline{p}\neq\infty$.
%\end{lemm}
%\noindent{\bf Proof of Lemma \ref{lemm_alpha=1}.} Since $u$ is a
%weak solution we have that $|u|^{1-\gamma}\in
%L^\infty(L^{2/(1-\gamma)})\cap L^{2/(1-\gamma)}(L^{6/(1-\gamma)})$.
%Interpolation gives that $|u|^{1-\gamma}\in L^{r}(L^s)$ for any
%$r\geq 2/(1-\gamma)$ such that:
%$$
%\frac{2}{r}+\frac{3}{s}=\frac{3(1-\gamma)}{2}.
%$$
%So: $||u|\Div(u/|u|)|\in L^{\overline{p}}$ for:
%$$
%\frac{2}{\overline{p}}+\frac{3}{\overline{q}}=\frac{2}{p}+\frac{3}{q}+\frac{3(1-\gamma)}{2},
%$$
%and:
%$$
%\frac{1}{\overline{p}}\leq \frac{1-\gamma}{2}+\frac{1}{p}.
%$$
%
\begin{lemm}\label{lemm_pq}
  For every $r$, $2\leq r<
6$, there exists a constant $C$ such that for every $\beta>0$, and
every function $f$ lying in $L^2(\R^3)$ and such that $\nabla f$
lies in $L^2(\R^3)$, we have:
$$
\beta\|f\|^2_{L^r(\R^3)}\leq \frac{1}{4} \|\nabla
f\|^{2}_{L^2(\R^3)}+C\beta^{\frac{1}{\theta}}\|f\|^{2}_{L^2(\R^3)},
$$
for $\theta=3/r-1/2$.
\end{lemm}
\noindent{\bf Proof of Lemma \ref{lemm_pq}.} Sobolev inequality
gives:
$$
\|f\|_{L^6(\R^3)}\leq C\|\nabla f\|_{L^2(\R^3)}.
$$
Interpolation gives:
$$
\beta\|f\|^2_{L^r(\R^3)}\leq
\left(\beta^{1/\theta}\|f\|^{2}_{L^2(\R^3)}\right)^\theta\left(\|f\|^2_{L^6(\R^3)}\right)^{1-\theta},
$$
where:
$$
\frac{\theta}{2}+\frac{1-\theta}{6}=\frac{1}{r},
$$
that is $\theta=3/r-1/2$. We end the proof using Minkowski
inequality:
$$
ab\leq \frac{a^p}{p}+\frac{b^q}{q},
$$
with:
$$
\theta=1/p\qquad 1-\theta=1/q,
$$
and:
$$
a=\frac{\left(\beta^{1/\theta}\|f\|^2_{L^2}\right)^\theta}{\varepsilon}\qquad
b=\varepsilon\left(\|\nabla f\|^2_{L^2}\right)^{1-\theta},
$$
for $\varepsilon$ small enough.
 \qquad\qed

We consider now $u$, a Leray-Hopf solution to Navier-Stokes
equation. Since  $u$ lies in $L^{2}(0,\infty; L^6(\R^3))$, for
almost every $t_0>0$, $u(t_0,\cdot)$ lies in $L^6(\R^3)$. It is
classical that it exists $T>t_0$ such that $u$ is smooth on
$(t_0,T)\times\R^3$. Moreover, from the Serrin's criterion, if
$T<\infty$, then:
$$
\lim_{t\to T}\|u\|_{L^3(t_0,t;L^9(\R^3))}=\infty.
$$
We will show that it cannot be the case. Note that $u(t_0,\cdot)\in
L^2(\R^3)\cap L^6(\R^3)$, so it lies in $L^3(\R^3)$. We consider $u$
on $(t_0,T)\times\R^3$.
% \noindent {\bf Proof of Theorem \ref{th_rough_theorem}.}
Multiplying  (\ref{eq_NS}) by $u|u|$, and integrating in $x$ we
find:
\begin{eqnarray*}
&&\frac{d}{dt}\int_{\R^3}\frac{|u |^3}{3}\,dx+\int_{\R^3}|u|(|\nabla
u|^2+|\nabla|u||^2)\,dx-\int_{\R^3}Pu\cdot\nabla|u|\,dx=0.
\end{eqnarray*}
Noting that:
$$
-\Delta P=\sum_{ij}\partial_i\partial_j(u_iu_j),
$$
we have for every $4/3<r<\infty$:
$$
\|P\|_{L^{3r/4}(\R^3)}\leq C_r\|u\|^2_{L^{3r/2}(\R^3)}.
$$
Since $\Div (u/|u|)\in L^p(L^q)$ for $2/p+3/q\leq1/2$, $q\geq 6$,
and $u\in L^a(L^b)$ for $2/a+3/b=3/2$, $2\leq b\leq6$, there exists
$\overline{p}>1$ and $2<\overline{q}<6$, such that
$|u|\Div(u/|u|)\in L^{\overline{p}}(L^{\overline{q}})$ with:
$$
\frac{1}{\overline{p}}=\frac{1}{p}+\frac{1}{a}\qquad
\frac{1}{\overline{q}}=\frac{1}{q}+\frac{1}{b}.
$$
Note that $2\leq \overline{q}<6$ and:
\begin{equation}\label{eq_2}
\frac{2}{\overline{p}}+\frac{3}{\overline{q}}\leq2.
\end{equation}
So, using (\ref{eq_incompressibility}), we have for every fixed time
$t$:
\begin{eqnarray*}
&&\frac{d}{dt}\int_{\R^3}\frac{|u|^3}{3}\,dx+\int_{\R^3}|u||\nabla|u||^2\,dx\leq
\int_{\R^3}|P||u|\left|\frac{u}{|u|}\cdot\nabla|u|\right|\,dx\\
&&\qquad\leq C\|u\|^3_{L^{3r/2}(\R^3)}\||u|\Div
(u/|u|)\|_{L^{\overline{q}}(\R^3)},
\end{eqnarray*}
with:
$$
\frac{2}{r}+\frac{1}{\overline{q}}=1.
$$
Using Lemma \ref{lemm_pq} with:
$$
f=|u|^{3/2},\qquad \nabla f=\frac{3}{2}|u|^{1/2}\nabla |u|,
$$
We find that:
\begin{eqnarray*}
&&C\|u\|^3_{L^{3r/2}(\R^3)}\||u|\Div
(u/|u|)\|_{L^{\overline{q}}(\R^3)}\\
&&\qquad=C\|f\|^2_{L^{r}(\R^3)}\||u|\Div
(u/|u|)\|_{L^{\overline{q}}(\R^3)}\\
&&\qquad\leq \frac{9}{16}\||u|^{1/2}\nabla
|u|\|^2_{L^2(\R^3)}+C\||u|\Div
(u/|u|)\|^{1/\theta}_{L^{\overline{q}}(\R^3)}\|u\|^3_{L^3(\R^3)},
\end{eqnarray*}
where:
$$
\theta=\frac{3}{r}-\frac{1}{2}=\frac{1}{2}\left(2-\frac{3}{\overline{q}}\right).
$$
From (\ref{eq_2}), this gives $1/\theta\leq \overline{p}$, hence
$\||u|\Div (u/|u|)\|^{1/\theta}_{L^{\overline{q}}(\R^3)}$ lies in
$L^1(0,T)$ with:
$$
\frac{d}{dt}\int_{\R^3}\frac{|u|^3}{3}\,dx+\frac{7}{16}\int_{\R^3}|u||\nabla|u||^2\,dx\leq
C \||u|\Div
(u/|u|)\|^{1/\theta}_{L^{\overline{q}}(\R^3)}\int_{\R^3}\frac{|u|^3}{3}\,dx.
$$
Gronwall argument gives that
$$
\lim_{t\to T}\int_{\R^3}|u|^3\,dx<\infty,
$$
and so:
\begin{eqnarray*}
&&\int_{t_0}^T\int_{\R^3}|u||\nabla|u||^2\,dx\,dt\\
&&\qquad=\frac{4}{9}\int_{t_0}^T\int_{\R^3}|\nabla|u|^{3/2}|^2\,dx\,dt,
\end{eqnarray*}
is finite too. Sobolev imbedding gives that $u\in
L^3(t_0,T;L^9(\R^3))$ which gives the desired contradiction. This
shows that $u$ is smooth on $(t_0,\infty)\times\R^3$ for almost
every $t_0>0$. The result of Theorem \ref{th_rough_theorem}
follows.\qquad\qed

\bibliography{biblio}

\end{document}